\documentclass[12pt,a4wide,twoside]{report}

\usepackage{amsthm,amssymb,amsmath}
\usepackage{mathrsfs,setspace,pstcol}
\usepackage{IEEEtrantools}
\usepackage{play}
\usepackage{epsfig}
\usepackage[nottoc]{tocbibind}

\input xy
\xyoption{all}

\theoremstyle{plain}
\newtheorem{theorem}{Theorem}[section]
\newtheorem{lemma}[theorem]{Lemma}
\newtheorem{corollary}[theorem]{Corollary}

\theoremstyle{definition}
\newtheorem{definition}[theorem]{Definition}
\newtheorem{example}[theorem]{Example}

\theoremstyle{remark}
\newtheorem{remark}[theorem]{Remark}

\newcommand{\spmk}{$SP_m(k)$}
\begin{document}


\begin{titlepage}
\enlargethispage{2cm}

\begin{center}

\vspace*{-1cm}

\textbf{\Large Topics in Ramsey Theory}\\[10pt]

\vspace*{2cm}

A Project Report Submitted \\
in Partial Fulfilment of the Requirements  \\
for the Degree of  \\
\vspace{5mm}
{\Large \bf MASTER OF SCIENCE } \\
in \\
{\large \bf Mathematics } \\

\vspace{10mm}
{\em  by} \\ \vspace{3mm}
{\large \bf Mano Vikash J} \\
{\large \bf (Roll No. IMS09036)}\\[.35in]

\vfill

\begin{figure}[h]
  \begin{center}
    \includegraphics[width=0.3\textwidth]{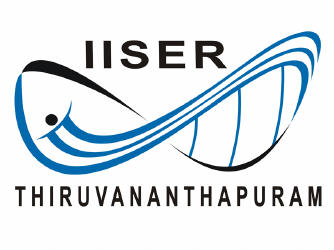}
  \end{center}
\end{figure}
\vspace*{0.25cm}

{\em\large to }\\
{\bf\large SCHOOL OF MATHEMATICS} \\
 {\bf\large {INDIAN INSTITUTE OF SCIENCE EDUCATION AND RESEARCH }}\\
{\bf\large THIRUVANANTHAPURAM - 695 016, INDIA}\\
{\it\large April 2014}
\end{center}

\end{titlepage}

\clearpage


\begin{center}

{\large{\bf{DECLARATION}}}

\end{center}

\pagenumbering{roman} \setcounter{page}{1}

\noindent

I declare that the matter embodied in this report :\textbf{``Topics in Ramsey Theory''} is the result of investigations carried out by me in the School of Mathematics, Indian Institute of Science Education and Research, Thiruvananthapuram, India under the supervision of Dr. Sujith Vijay. It has not been submitted elsewhere
for the award of any degree. In keeping with the general practice of reporting scientific observations, due acknowledgement has been made whenever the work described is based on the findings of other investigations. Any omission which might have occured by oversight or error in judgement is regretted.

\vspace{4cm}

\noindent Thiruvananthapuram - 695 016 \hfill (Mano Vikash J)

\noindent April  2014  \hfill (IMS 09036)

\clearpage

\pagenumbering{roman} \setcounter{page}{2}
\begin{center}
{\large{\bf{CERTIFICATE}}}
\end{center}

\noindent
This is to certify that the work contained in this project report
entitled \textbf{``Topics in Ramsey Theory''}  submitted
by \textbf{Mano Vikash J} (\textbf{Roll No: IMS09036}) to Indian Institute of Science Education and Research Thiruvananthapuram
towards partial requirement of {\bf Master of Science} in Mathematics   has been carried out
by him under my supervision and that it has not been submitted elsewhere
for the award of any degree.

\vspace{4cm}

\noindent Thiruvananthapuram - 695 016 \hfill (Dr. Sujith Vijay)

\noindent April  2014 \hfill Project Supervisor

\clearpage

\begin{center}
{\Large{\bf{Abstract}}}
\end{center}

Ramsey theory is the study of conditions under which mathematical objects show order when partitioned. Ramsey theory on the integers concerns itself with partitions of $[1,n]$ into $r$ subsets and asks the question whether one (or more) of these $r$ subsets contains a $k$-term member of $\mathcal{F}$, where $[1,n]=\{1,2,3,\ldots,n\}$ and $\mathcal{F}$ is a certain family of subsets of $\mathbb{Z}^+$. When $\mathcal{F}$ is fixed to be the set of arithmetic progressions, the corresponding Ramsey-type numbers are called the van der Waerden numbers. 

I started the project choosing $\mathcal{F}$ to be the set of semi-progressions of scope $m$. A semi-progression of scope $m\in \mathbb{Z}^+$ is a set of integers $\{x_1,x_2,\ldots,x_k\}$ such that for some $d\in\mathbb{Z}^+$, $x_{i}-x_{i-1}\in\{d,2d,\ldots,md\}$ for all $i\in\{2,3,\ldots,k\}$. The existence of Ramsey-type numbers corresponding to semi-progressions follows immediately from the existence of van der Waerden numbers. However, their exact values are not known. We use $SP_m(k)$ to denote these numbers as a Ramsey-type function of $k$ for a fixed scope $m$. The best known lower bound for this function $SP_m(k)$ was a second degree polynomial. During this project, I used the probabilistic method to increase this to an exponential lower bound for any fixed $m$. That is $SP_m(k)>c^k$, for some $c>1$. The base of the exponential $c$ is a strictly decreasing function of $m$ that tends to $1$ as $m$ tends to infinity. The first chapter starts with a brief introduction to Ramsey theory and then explains the problem considered. In second chapter, I give the results obtained on semi-progressions.

In the third chapter, I will discuss the lower bound obtained on $Q_1(k)$. When $\mathcal{F}$ is chosen to be quasi-progressions of diameter $n$, the corresponding Ramsey-type numbers obtained are denoted as $Q_n(k)$. A quasi-progression of diameter $n\in \mathbb{Z}^+$ is a set of integers $\{x_1,x_2,\ldots,x_k\}$ such that for some $d\in\mathbb{Z}^+$, $x_{i}-x_{i-1}\in\{d,d+1,\ldots,d+m\}$ for all $i\in\{2,3,\ldots,k\}$. The approach only gives an exponential lower bound for $Q_1(k)$. This bound obtained beats the previous best known bound. The approach does not work when $n>1$ because we only get $Q_n(k)>c^k$ for some $c<1$.

The last chapter gives an exposition of advanced probabilistic techniques, in particular concentration inequalities. When using the probabilistic method, these inequalities will be extremely useful to estimate the concentration of a random variable. In this chapter, I will discuss some advanced concentration inequalities and how to apply them. I will mostly restrict myself to applications to graph theoretic problems.

\clearpage

\begin{center}
{\Large{\bf{Acknowledgements}}}
\end{center}

Firstly, I thank my supervisor, Dr. Sujith Vijay for his technical guidance, encouragement and time. His intuition and guidance have had a very significant impact on the work presented here. At many stages during this project, I benefited immensely from his advice. I am also grateful for his careful editing. This project has been not only an enjoyable experience for me, but also genuinely helpful in finding my research interests. I am extremely lucky to have had all the interesting discussions about academics with him during the course of the project. I also express my gratitude for the courses he taught in IISER Thiruvananthapuram, which kindled my interest. It was a privilege to learn from him. I extend my gratitude to all the faculty in IISER Thiruvananthapuram, especially the ones in the School of Mathematics, for their support and the courses they taught me. I thoroughly enjoyed all of them. I thank my family and friends for their support and encouragement.

\tableofcontents
\clearpage

\newpage

\pagenumbering{arabic}
\setcounter{page}{1}


\chapter{Introduction}

Ramsey theory, named after British mathematician Frank Plumpton Ramsey, studies preservation of properties of mathematical object under set partitions. Although Ramsey's theorem itself dealt with study of properties of graphs, there were earlier results of similar flavor proved on integers. Ramsey theory on the integers studies preservation of properties of of set of integers when partitioned. As this report is primarily concerned with results in integer Ramsey theory, we restrict ourselves to this right from the beginning.

\section{Introduction to Ramsey Theory on the Integers}

Let us start with some definitions.

\subsection{Preliminaries}
\begin{definition}\label{ap}
A $k$-term arithmetic progression is a sequence of numbers $x_1,x_2,\ldots,x_k$ where $x_{i}-x_{i-1}=d$ for $2\leq i\leq k$
\end{definition}
We now look at a generalisation of arithmetic progressions.
\begin{definition}\label{sp}
A $k$-term semi-progression of scope $m$ is a sequence of numbers $x_1,x_2,\ldots,x_k$ such that for some $d\in \mathbb{Z}^+$, $x_{i}-x_{i-1}\in \{d, 2d,\ldots,md\}$ for $2\leq i\leq k$.
\end{definition}

\begin{remark}
Note when $m=1$, the definition of semi-progression is the same as that of an arithmetic progression. For any $m>1$, the set of semi-progressions of scope $m$ contains the set of arithmetic progressions. But every semi-progression need not be an arithmetic progression. For example, though $1,2,3,5,6,7$ forms a semi-progression of scope $2$, it is not an arithmetic progression.
\end{remark}

The second chapter, I will only consider semi-progressions. But in the third chapter, we will be considering quasi-progressions too. With that in mind, I make the following definition.
\begin{definition}\label{sp}
A $k$-term quasi-progression of diameter $n$ is a sequence of numbers $x_1,x_2,\ldots,x_k$ such that for some $d\in \mathbb{Z}^+$, $x_{i}-x_{i-1}\in \{d, d+1,\ldots,d+n\}$ for $2\leq i\leq k$.
\end{definition}

Let $[1,n]=\{1,2,\ldots,n\}$. Consider a partition of $[1,n]$ into $r$ disjoint subsets. This partitioning can be equivalently thought of as a coloring of $[1,n]$ where each subset is assigned a unique color among $r$ possible colors. The process of coloring can be formalised as follows.

\begin{definition}\label{coloring}
An $r$-coloring of $[1,n]$ is a function $\chi: [1,n]\rightarrow C$ where $|C|=r$.
\end{definition}
$C$ itself is usually taken to be $\{1,2,\ldots,r\}$. If we consider the elements of $[1,n]$ mapped to $i$ to be $S_i$, the $r$ subsets $S_1,S_2,\dots,S_r$ will give back our partition of $[1,n]$. 
\begin{definition}\label{monochromatic}
Given any $G\subset [1,n]$, we say that $\chi$ is monochromatic on the set $G$ if $\chi$ is constant on G.
\end{definition}

\subsection{Van der Waerden's numbers}
For a given $r\in\mathbb{Z}^+$, there are $r^n$ possible colorings of $[1,n]$. We fix a $k\in\mathbb{Z}^+$ and ask the question whether there exists an $n$ for which we can say that each of the $r^n$ colorings has a monochromatic $k$-term arithmetic progression.\footnote{Note that once we fix $r$ and $k$, if the result is true for $n$, the result will be true for any $m>n$.}
It turns out that there exists such an $n$ for any $r,k\in \mathbb{Z}^+$. This is one of the fundamental results of Ramsey theory and was proved by van der Waerden. I state his result here; its proof can be found in \cite{landman2004ramsey}.
\begin{theorem}
  \label{thm:vanderwaerden}
  Let $k,r\geq2$ be integers. There exists a least positive integer $w(k;r)$ such that for every $r$-coloring of $[1,w(k;r)]$, there is a monochromatic arithmetic progression of length $k$.
\end{theorem}

Once the existence of $w(k;r)$ is established, the next step is to find their value. But this has not been easy. Only very few nontrivial van der Waerden numbers have been found so far. This is because, in order to find a van der Waerden number exactly, we would have to check $r^{w(k;r)}$ colorings. This function grows very fast and it will take thousands of years for any present day computer to perform such a check (except for a few possible cases of small $k$ and $r$). Given that van der Waerden numbers are difficult to find exactly, the next step is to try to find bounds on these numbers. To show that $f(k;r)$ is a lower bound for $w(k;r)$, we only need to show that some $r$-coloring of $[1,f(k;r)]$ does not have any monochromatic $k$-term arithmetic progression. Some of the best known lower bounds are

\begin{theorem}
  \label{lb1}
  Let $p\geq5$ and $q $ be primes. Then
  \begin{equation*}
  w(p+1;q)\geq p(q^p-1)+1
  \end{equation*}
\end{theorem}

\begin{theorem}
  \label{lb2}
  For all $r\geq2$,
  \begin{equation*}
  w(k;r)> \frac{r^k}{ekr}(1+o(1))
  \end{equation*}
\end{theorem}

Upper bounds are more difficult to find. This is because in order to show that $g(k;r)$ is an upper bound for $w(k;r)$, we have to prove that every $r$-coloring of $[1,g(k;r)]$ has a monochromatic $k$-term arithmetic progression. The best known upper bound, proved by Gowers in \cite{gowers}, is stated below

\begin{theorem}
  \label{ub}
  For $k,r\geq2$,
  \begin{equation*}
  w(k;r)\leq 2^{2^{r^{2^{2^{k+9}}}}}
  \end{equation*}
\end{theorem}

\section{The semi-progression problem}

I started the project studying Ramsey-type functions corresponding to semi-progressions and tried to find lower bounds for them. This section describes the Ramsey-type function studied and also states the main results known about them.
\subsection{Ramsey-type function for semi-progression ($SP_m(k)$)}
In the previous section, after the introduction to Ramsey theory on the integers, we looked at a specific example of this, namely the van der Waerden numbers ($w(k;r)$). In this section, we will look at a generalisation of this. Van der Waerden number is the minimum number $w(k;r)$ such that any $r$-coloring of $[1,w(k;r)]$ contains a monochromatic $k$-term arithmetic progression. Now, why should one restrict to arithmetic progressions? We can generalise the definition to semi-progressions and quasi-progressions as follows:
\begin{definition}
  \label{sp}
  Let $k\geq2$ and $m\geq1$ be positive integers. Then, $SP_m(k)$ is the minimum positive integer such that any $2$-coloring of $[1,SP_m(k)]$ will contain a monochromatic $k$-term semi-progression of scope $m$.
\end{definition}

\begin{definition}
  \label{qp}
  Let $k\geq2$ and $n\geq0$ be positive integers. Then, $Q_n(k)$ is the minimum positive integer such that any $2$-coloring of $[1,Q_n(k)]$ will contain a monochromatic $k$-term quasi-progression of diameter $n$.
\end{definition}

Let us restrict ourselves to $SP_m(k)$ for now.

\begin{remark}
  The existence of $SP_m(k)$ follows directly from the existence of van der Waerden numbers. Because every arithmetic progression is a semi-progression and there exist $w(k;2)$ such that any $2$-coloring of $[1,w(k;2)]$ will contain a monochromatic $k$-term arithmetic progression, this $[1,w(k,2)]$ will automatically contain a $k$-term semi-progression for any fixed scope $m\in \mathbb{Z}^+$. Hence, we can say that 
  \begin{equation*}
  SP_m(k)\leq w(k;2)
  \end{equation*}
  
  Further, for $m_1,m_2\in \mathbb{Z}^+$ and $m_1>m_2$, any semi-progression of scope $m_2$ is also a semi-progression of scope $m_1$. This would imply that
  \begin{equation*}
  w(k;2)= SP_1(k)\geq SP_2(k)\geq SP_3(k)\geq \ldots
  \end{equation*}
\end{remark}
 
\subsection{Some known results on \spmk}
I end this chapter by stating some of the best bounds known for \spmk. They were proved by Landman in \cite{bookreference}. The upper bounds for $SP_m(k)$ are easier to find compared to finding upper bounds for $w(k;2)$. The following theorem states the upper bound under some restrictions on $m$ and $k$.
\begin{theorem}
  \label{ubsp}
  Let $m\geq2$. Assume $m<k<2m$. Let $c=\left\lceil {\frac{m}{2m-k}}\right\rceil $. Then,
  \begin{equation*}
  SP_m(k)\leq 2c(k-1)+1
  \end{equation*}
\end{theorem}

Next, we turn to lower bounds. The following theorem gives the best constructive result.
\begin{theorem}
  \label{lbsp}
  Let $k\geq2$ and $m\geq1$. Let $\lambda(k,m)=\left\lceil{\frac{k-1}{\lceil k/m\rceil } }\right\rceil $. Then,
  \begin{equation*}
  SP_m(k)\geq 2(k-1)\left(\left\lceil\frac{k}{\lambda(k,m)} \right\rceil -1\right)+1
  \end{equation*}
\end{theorem}
The next chapter gives a probabilistic lower bound obtained during this project.

\chapter{Lower bounds obtained for $SP_m(k)$}

The probabilistic method was used to get exponential lower bounds for \spmk. I will give a brief introduction to probabilistic method before stating the results obtained during the project.

\section{The Probabilistic Method}
The probabilistic method was pioneered by Paul Erd\H{o}s. This method can be used to prove the existence of certain mathematical objects. I will illustrate the method by proving a simple result for the van der Waerden numbers found in \cite{alon2011probabilistic,jukna2011extremal}. The cited books contain more examples of application of the probabilistic method.
\begin{theorem}
\begin{equation*}
w(k;2) \geq  \left( \frac{2^kk}{2}\right)^{1/2}
\end{equation*}
\end{theorem}

\begin{proof}
Fix $n\in \mathbb{Z}^+$. Let us color the elements of $[1,n]$ randomly red or blue with probability $1/2$ and call this coloring $\chi$. Now, fix a particular arithmetic progression with $k$ terms (and call this arithmetic progression $i$). Let $A_i$ be the event that $i$ is monochromatic under this coloring. Then, the probability of the event $A_i$ is $2^{1-k}$. The previous statement can be written as $P[A_i]=2^{1-k}$.

Now, let us try to get an upper bound on the number of arithmetic progressions in the set $\{1,2,\ldots,n\}$. These progressions are specified uniquely by its initial value and common difference. There are at most $n$ choices for the initial value and at most $n/k$ choices for the common difference. Hence, there are at most $n^2/k$ ways to choose $i$. Thus, the probability that at least one of them is monochromatic is given by
\begin{IEEEeqnarray*}{rCl}
P[\vee_i {A}_i]&\leq& \sum_i P[A_i]\\ &\leq & \frac{n^2}{k}\times2^{1-k} \IEEEyesnumber \label{eqn1}
\end{IEEEeqnarray*}
If the right hand side of equation \ref{eqn1} is less than 1, then the probability that none of the events $A_i$ occurs is non-zero ($P[\wedge_i \overline{A}_i]>0$). For this to happen,  
\begin{IEEEeqnarray*}{rCl}
\frac{n^2}{k}\times2^{1-k}&<&1\\ \Rightarrow n & < & \left( \frac{2^kk}{2}\right)^{1/2}
\end{IEEEeqnarray*}
So, when $n < \left( \frac{2^kk}{2}\right)^{1/2}$, there exist colorings of $[1,n]$ without any $k$-element arithmetic progression being monochromatic. Thus,
\begin{equation}
w(k;2)\geq  \left( \frac{2^kk}{2}\right)^{1/2}
\end{equation}
\end{proof}
Note that the probabilistic method only gives an existential result. In this case, it gave a lower bound for $w(k;2)$ without giving an explicit coloring for $[1,w(k;2)]$ that avoids a monochromatic $k$-term arithmetic progression.

\section{Lower bounds obtained}
This section gives the main results obtained during the project. I will prove the result for scope $m=2$ and then extend the proof for any scope $m$. I will also extend the approach to $Q_1(k)$ in the next chapter.

\subsection{Scope $2$}

\begin{theorem}
\begin{equation*}
\sqrt{\frac{3k}{4}}\left(\sqrt{\frac{4}{3}}\right)^k\leq SP_2(k)
\end{equation*}
\end{theorem}

\begin{proof}
Let $N$ be an integer (to be picked later). Let $[N]=\{1,2,\ldots,N\}$. We pick the elements of $[N]$ one by one and color them independently red or blue. The total number of possible colorings is $2^N$. Now, let us count the number of colorings that have a monochromatic $k$-term semi-progression of scope 2. Let us call this number $S$.
For a fixed $a$ ($a=x_1$ is the first term of the semi-progression) and $d$, let $T_{a,d}$ denote the number of colorings of $[N]$ that have a monochromatic semi-progression of scope 2 with first term $a$ and difference $d$. 
\begin{equation*}
S\leq \sum_{a,d}T_{a,d}
\end{equation*}
We have atmost $N-k+1$ possible choice for $a$ and atmost $N/k$ possible choice for $d$. Hence, 
\begin{equation*}
S\leq \sum_{a,d}T_{a,d}\leq \frac{(N-k+1)N}{k} \max_{a,d}T_{a,d}
\end{equation*}
Now, let us try to find $\max_{a,d}T_{a,d}$ (call this quantity $\lambda$). First note that $T_{a,d}$ is not same for all choices of $a$ and $d$. Suppose $a=N-k+1$ and $d=1$, then we have only one possible semi-progression with this particular $a$ and $d$. If this is monochromatic, we have $N-k$ other elements that can be coloured any other way. This progression itself can take any of two colours. Hence, $T_{a,d}=2^{N-k+1}$. But, for other choices of $a$ and $d$, we may have semi-progressions that take one or more jumps of $2d$. 

Given a semi-progression of scope $2$ and length $k$, it takes $k-1$ jumps and these jumps may be of size $d$ or $2d$. Let $\gamma_{a,d}$ be the set of all colorings that contains a monochromatic semi-progression of scope $2$ with starting term $a$ and difference $d$. $\lambda=\max_{a,d}|\gamma_{a,d}|$. Let
\begin{equation*}
R_{a,d}=\{a+qd\quad:\quad q\in\{0,1,2,3,\ldots,(k-1)+r\}\}
\end{equation*}
where $0\leq r\leq k-1$ and $r$ is the maximum number in this range such that $a+(k-1+r)d\leq N$. Note that $r$ will depend on $a$ and $d$. If for some $a$ and $d$, $a+(k-1)d>N$, we will not be able to get a $k$-term semi progression with this $a$ and $d$ that does not exceed $N$.  Also note that if we pick any coloring in $\gamma_{a,d}$ and change the color of any element in $[N]\setminus R_{a,d}$, the resulting coloring will still be in $\gamma_{a,d}$. There are $2^{N-(k-1)-r-1}$ ways to color the elements in $[N]\setminus R_{a,d}$. Let us count modulo these redundant elements and add this factor finally to our result. Now, we only need to worry about the color of the elements of $R_{a,d}$ in each of the colorings in $\gamma_{a,d}$.

Given any coloring in $\gamma_{a,d}$, we look at the least number of jumps of size $2d$ we need to take among the $k-1$ jumps starting from $a$ to get a monochromatic $k$-term semi progression. Let this number be $l$. We map this coloring to this number. We now count the number of colorings mapped to each number and sum them up. For a given number $l$, the number of colorings mapped to $l$ is given by ${k-1\choose l}2^{r+1-l}$. Hence, we get that 
\begin{IEEEeqnarray*}{rCl}
\frac{|\gamma_{a,d}|}{2^{N-(k-1)-r-1}}&=& {k-1\choose 0}2^{r+1}+ {k-1\choose 1}2^{r}+\ldots+ {k-1\choose r}2^{1}\\
&\leq& 2^{r+1} \left(1+\frac{1}{2}\right)^{k-1}\\
&=&(2)^{r+1}\left(\frac{3}{2}\right)^{k-1}
\end{IEEEeqnarray*}
Hence,
\begin{IEEEeqnarray*}{rCl}
{|\gamma_{a,d}|}&\leq& {2^{N-(k-1)-r-1}} (2)^{r+1}\left(\frac{3}{2}\right)^{k-1}\\
&=&2^{N-2k+2}3^{k-1}
\end{IEEEeqnarray*}
Hence,
\begin{IEEEeqnarray*}{rCcCl}
\lambda&=&\max_{a,d}|\gamma_{a,d}|&\leq&2^{N-2k+2}3^{k-1}
\end{IEEEeqnarray*}
Hence, we get
\begin{IEEEeqnarray*}{rCl}
S&\leq&(N-k+1)\frac{N}{k} 2^{N-2k+2}3^{k-1}\\
&\leq&4 \frac{N^2}{k}2^N\frac{3^k}{2^{2k}}\frac{1}{3}
\end{IEEEeqnarray*}
If this quantity is less than the total number of colorings ($2^N$), then we can say that there exist colorings of $[N]$ that do not have any monochromatic $k$-term semi-progressions of scope $2$. Thus, we impose this as a condition and get some bound for $N$.
\begin{IEEEeqnarray*}{rCl}
4 \frac{N^2}{k}2^N\frac{3^k}{2^{2k}}\frac{1}{3}&<&2^N\\
N^2&<&\frac{3k}{4}\frac{2^{2k}}{3^k}\\
N^2&<&\frac{3k}{4}\left(\frac{4}{3}\right)^k\\
N&<&\sqrt{\frac{3k}{4}}\left(\sqrt{\frac{4}{3}}\right)^k
\end{IEEEeqnarray*}
When $N$ satisfies the above inequality, we have colorings of $[N]$ that do not have any monochromatic $k$-term semi-progression of scope $2$.
\end{proof}

\subsection{Arbitrary scope $m\geq 2$}
\begin{theorem}
\begin{equation*}
\sqrt{\frac{(2^m-1)k}{2^m}}\left(\sqrt{\frac{2^m}{2^m-1}}\right)^k\leq SP_m(k)
\end{equation*}
\end{theorem}
\begin{proof}
Let $N$ be an integer. We pick the elements of $[N]$ one by one and color them independently red or blue. Let $S$ denote the number of colorings that have a monochromatic $k$-term semi-progression of scope $m$. For a fixed $a$ and $d$, let $T_{a,d}$ denote the number of colorings of $[N]$ that have a monochromatic semi-progression of scope $m$ with first term $a$ and difference $d$. Then,
\begin{equation*}
S\leq \sum_{a,d}T_{a,d}
\end{equation*}
We have atmost $N-k+1$ possible choice for $a$ and atmost $N/k$ possible choice for $d$. Hence, 
\begin{equation*}
S\leq \sum_{a,d}T_{a,d}\leq \frac{(N-k+1)N}{k} \max_{a,d}T_{a,d}
\end{equation*}
Let $\max_{a,d}T_{a,d}=\lambda$. Let $\gamma_{a,d}$ be the set of all colorings that contains a monochromatic semi-progression of scope $m$ with starting term $a$ and difference $d$. Then, $\lambda=\max_{a,d}|\gamma_{a,d}|$. Let
\begin{equation*}
R_{a,d}=\{a+qd\quad:\quad q\in\{0,1,2,3,\ldots,(k-1)+r\}\}
\end{equation*}
where $ 0\leq r\leq (m-1)(k-1)$.  $r$ is the maximum number in this range such that $a+(k-1+r)d\leq N$. Note that $r$ depends on $a$ and $d$. If for some $a$ and $d$, $a+(k-1)d>N$, we will not be able to get a $k$-term semi progression with this $a$ and $d$ that does not exceed $N$. Also note that if we pick any coloring in $\gamma_{a,d}$ and change the color of any element in $[N]\setminus R_{a,d}$, the resulting coloring will still be in $\gamma_{a,d}$. There are $2^{N-(k-1)-r-1}$ ways to color the elements in $[N]\setminus R_{a,d}$. Let us count modulo these redundant elements. Now, we only need to worry about the color of the elements of $R_{a,d}$ in each of the colorings in $\gamma_{a,d}$.

Given a coloring in $\gamma_{a,d}$, there exists $(x_1,x_2,\ldots,x_m)$ such that \\\mbox{$x_1+x_2+\ldots+x_m=k-1$} and we would be able to take $x_i$ jumps of size $id$ (for each $i$) starting from $a$ and get a monochromatic semi-progression in this coloring. Among all such $m$-tuples, pick the lexicographically highest one and map this coloring to this $m$-tuple. Now, each coloring in $\gamma_{a,d}$ is mapped to a $m$-tuple. All we need to do now is to count the number of colorings mapped to each $m$-tuple and sum them up. For a given $m$-tuple, the number of colorings mapped to this (modulo the redundant elements) is
\begin{equation*}
{k-1\choose x_1,x_2,\ldots,x_m}2^{r+1-\sum_{i=1}^m(i-1)x_i}
\end{equation*}
Hence,
\begin{IEEEeqnarray*}{rCl}
\frac{|\gamma_{a,d}|}{2^{N-(k-1)-r-1}}&=& \sum_{\substack{(x_1,x_2,\ldots,x_m)\\ x_1+x_2+\ldots +x_m=k-1 \\ \sum_i (i-1)x_i\leq r}} {k-1\choose x_1,x_2,\ldots,x_m}2^{r+1-\sum_{i=1}^m (i-1)x_i}\\
&\leq& \sum_{\substack{(x_1,x_2,\ldots,x_m)\\ x_1+x_2+\ldots +x_m=k-1}} {k-1\choose x_1,x_2,\ldots,x_m}2^{r+1-\sum_{i=1}^m (i-1)x_i}\\
&=& 2^{r+1} \left(1+\frac{1}{2}+\frac{1}{2^2}+\ldots+\frac{1}{2^{m-1}}\right)^{k-1}\\
&=& 2^{r+1}\left(\frac{1-\left(\frac{1}{2}\right)^m}{1-\frac{1}{2}}\right)^{k-1}\\
&=& 2^{r+1} 2^{k-1}\left(\frac{2^m-1}{2^m}\right)^{k-1}
\end{IEEEeqnarray*}
\begin{IEEEeqnarray*}{rCl}
{|\gamma_{a,d}|}&\leq& {2^{N-(k-1)-r-1}} 2^{r+1} 2^{k-1}\left(\frac{2^m-1}{2^m}\right)^{k-1}\\
&=& 2^N \left(\frac{2^m-1}{2^m}\right)^{k-1}
\end{IEEEeqnarray*}
Hence,
\begin{IEEEeqnarray*}{rCcCl}
\lambda&=&\max_{a,d}|\gamma_{a,d}|&\leq& 2^N \left(\frac{2^m-1}{2^m}\right)^{k-1}
\end{IEEEeqnarray*}
Hence,
\begin{IEEEeqnarray*}{rCl}
S&\leq&(N-k+1)\frac{N}{k} 2^N \left(\frac{2^m-1}{2^m}\right)^{k-1} \\
&\leq& \frac{N^2}{k} 2^N \left(\frac{2^m-1}{2^m}\right)^{k-1}
\end{IEEEeqnarray*}
If this quantity is less than the total number of colorings ($2^N$), then we can say that there exist colorings of $[N]$ that do not have any $k$-term semi-progressions of scope $m$. Thus, we impose this condition and get some bound for $N$.
\begin{IEEEeqnarray*}{rCl}
\frac{N^2}{k} 2^N \left(\frac{2^m-1}{2^m}\right)^{k-1}&<&2^N\\
N^2&<&{k}\left(\frac{2^{m}}{2^m-1}\right)^{k-1}\\
N&<&\sqrt{k} \left(\sqrt{\frac{2^m}{2^m-1}}\right)^{k-1}
\end{IEEEeqnarray*}
When $N$ satisfies the above inequality, we have colorings of $[N]$ that do not have any monochromatic $k$-term semi-progression of scope $m$.
\end{proof}

\chapter{Results for quasi-progressions}

Firstly, let us recall the following definitions. A quasi-progression of length $k$ and diameter $n$ is a sequence of integers $\{x_1,x_2,\ldots ,x_k\}$ such that $x_i-x_{i-1}\in\{d,d+1,\ldots,d+n\}$ for $i\in\{2,3,\ldots,k\}$. 
\begin{definition}
Let $k\geq1$ and $n\geq0$. $Q_n(k)$ is the least positive integer such that any 2-coloring of $[1,Q_n(k)]$ will have a monochromatic $k$-term quasi-progression of diameter $n$.
\end{definition}

\section{Best known bounds}
When the diameter is large, the values of $Q_n(k)$ are known exactly. The following result (proved in \cite{jobson}) gives the exact values for $Q_n(k)$ when the diameter is large.
\begin{theorem}
\begin{equation*}
Q_{k-i}(k)=2ik-4i+2r-1
\end{equation*}
if $k=mi+r$ for integers $m,r$ such that $3 \le r< \frac{i}{2}$ and $r-1 \le m$.
\end{theorem}
The best known lower bounds for $Q_n(k)$ are second degree polynomials when $n>1$. For $n=1$, we have the following result proved by Vijay in \cite{Vijay2010}.
\begin{theorem}
$$Q_1(k)>\beta^k$$
where $\beta$ is the smallest positive real root of the equation \\ $y^{24}+8y^{20}-112y^{16}-128y^{12}+1792y^8+1024y^4-4096=0$.
\end{theorem}
This root comes out to be $\beta=1.08226..$. Next, we turn to upper bounds. Landman \cite{raey} proved upper bounds for all $n\geq \left\lceil\frac{2k}{3}\right\rceil$. He proved that they are bounded by polynomials. In particular, we are interested in quasi-progressions of small diameter. For $n= \left\lceil\frac{2k}{3}\right\rceil$, he got the following result.
\begin{theorem}
\begin{equation*}
Q_{\left\lceil\frac{2k}{3}\right\rceil}(k)\leq \frac{43}{324}k^3(1+o(1))
\end{equation*}
\end{theorem}

\section{Results obtained}

Now, we prove an exponential lower bound for $Q_1(k)$. This marginally improves the current best known bound (see \cite{Vijay2010}). We can also use this approach for general $r$-colorings (see \cite{manosujith2014}).

\begin{theorem}
\begin{equation*}
Q_1(k)\ge c_0 (1.0823)^k
\end{equation*}
for some $c_0\in\mathbb{R}^+$.
\end{theorem}
\begin{proof}
Let $N$ be an integer (to be picked later). Let $[N]=\{1,2,\ldots,N\}$. We pick the elements of $[N]$ one by one and color them independently red or blue. The total number of possible colorings is $2^N$. Now, let us count the number of colorings that have a $k$-term quasi-progression of diameter $1$. Let us call this number $S_k$.

For a fixed $a$ ($a=x_1$ is the first term of the quasi-progression) and $d$, let $T_{a,d,k}$ denote the number of colorings of $[N]$ that have monochromatic quasi-progressions of diameter $1$ with first term $a$ and difference $d$. 
\begin{equation*}
S\leq \sum_{a,d}T_{a,d}
\end{equation*}
We have atmost $N-k+1$ possible choice for $a$ and atmost $N/k$ possible choice for $d$. Hence, 
\begin{equation}
\label{maineqn}
S_k\leq \sum_{a,d}T_{a,d,k}\leq \frac{(N-k+1)N}{k} \max_{a,d}T_{a,d}
\end{equation}
Now, let us try to find $\max_{a,d}T_{a,d}$ (call this quantity $\Omega_k$). First note that $T_{a,d,k}$ is not same for all choices of $a$ and $d$. Suppose $a=N-k+1$ and $d=1$, then we have only one possible quasi-progression with this particular $a$ and $d$. If this is monochromatic, we have $N-k$ other elements that can be coloured any other way. This progression itself can take any of two colours. Hence, $T_{a,d}=2^{N-k+1}$. But, for other choices of $a$ and $d$, we may have quasi-progressions that take jumps of size $d+1$. So, the number is of coloring increases if the possibility of taking jumps of size $d+1$ is more. Hence, the maximum occurs when $k-1$ jumps of size $d+1$ is possible.

Given a quasi-progression of diameter $1$ and length $k$, it takes $k-1$ jumps and these jumps may be of size $d$ or $d+1$. Let $\gamma_{a,d,k}$ be the set of all colorings of $[N]$ that contains a monochromatic quasi-progression of diameter $1$ with starting term $a$ and difference $d$. $\Omega_k=\max_{a,d}|\gamma_{a,d,k}|$. Let us also assume that $N\geq2k-1$. Let
\begin{equation*}
R_{a,d}=\cup_{i=0}^{k-1}P_{a,d,i}
\end{equation*}
where $P_{a,d,i}=\{a+id+q\quad:\quad q\in\{0,1,2,3,\ldots,i\}\}$ for $0\leq i\leq k-2$ and $P_{a,d,k-1}=\{a+(k-1)d+q\quad:\quad q\in\{0,1,2,3,\ldots,t\}\}$ where $0\leq t\leq k-1$ and $t$ is the maximum number in this range such that $a+(k-1)d+t\leq N$. Clearly, $t$ is a function of $a$ and $d$. Note that if $a+(k-1)d>N$, then we will not be able to pick such a $t$. This means that it is impossible to find a $k$-term quasi progression with this $a$ and $d$ that does not exceed $N$. Since $N\geq2k-1$, we will be able to take $k-1$ jumps of size $d+1$ starting from $1$ with difference $1$. Hence, $\max_{a,d}T_{a,d,k}=T_{1,1,k}$. Hence, we only need to count $T_{1,1,k}$ 

Also note that if we pick any coloring in $\gamma_{1,1,k}$ and change the color of any element in $[N]\setminus R_{1,1}$, the resulting coloring will still be in $\gamma_{1,1}$. Let $s=|R_{1,1}|$. Then, we can say that $s=\frac{k(k-1)}{2}+t-w$ where $w$ takes care of the repetitions in $P_{a,d,i}$. There are $2^{N-s}$ ways to color the elements in $[N]\setminus R_{1,1}$. Let us count modulo these redundant elements and add this factor finally to our result. 
\begin{IEEEeqnarray*}{rCl}
\Omega_k&=&2^{N-s}\Psi_k\\
&=&2^{N-s}2^{s-k}\lambda_k
\end{IEEEeqnarray*}
where $2^{s-k}\lambda_k=\Psi_k$. 

We now map each coloring in $\gamma_{1,1,k}$ to a $(k-1)$-tuple where each element of the $(k-1)$-tuple belongs to $\{0,1\}$. Firstly, we pick a coloring in $\gamma_{1,1,k}$ (say $\chi$). This coloring contains atleast one monochromatic $k$-term quasi-progression with starting term $1$ and difference $1$. Take one such quasi-progression. We look at the $k-1$ jumps it takes one by one starting from the first. If the first jump is of size $1$, the first element in the $(k-1)$-tuple is $0$. If the jump is of size $2$, the first element in the $(k-1)$-tuple is $1$. Now, we look at the second jump. This will determine the second element in the $(k-1)$-tuple. We do this for the all the $(k-1)$ terms, we will get a $(k-1)$-tuple. Then, we do this for all the possible monochromatic quasi-progressions with starting term $1$ and difference $1$ in $\chi$. We map the $\chi$ to the lexicographically least among the possible $(k-1)$-tuples.

Note that each $(k-1)$-tuple contributes to $\lambda_k$. The presence of $(\ldots,1,0,\ldots)$ in the $(k-1)$-tuple gives us additional information about the coloring and hence decreases the contribution of this $(k-1)$-tuple to $\lambda_k$ by a factor of $1/2$ for each 0 followed immediately by 1 in the $(k-1)$-tuple.
\begin{IEEEeqnarray*}{rCl}
\lambda_k&=&\lambda_{k,0}+\lambda_{k,1}
\end{IEEEeqnarray*}
where $\lambda_{k,0}$ gives the contribution of $(k-1)$-tuples that end with $0$ and $\lambda_{k,1}$ gives the contribution of $(k-1)$-tuples that end with $1$.

$k$-tuples can be got from $(k-1)$-tuples by adding a $0$ or $1$ to the right of the $(k-1)$-tuple. If we add a $1$ to the right of the $(k-1)$-tuple, we will not get any additional information. But, if we add a $0$ to the right of the $(k-1)$-tuple, we will get additional information about the coloring if the initial $(k-1)$-tuple ended with a $1$. Hence, we can write
\begin{equation*}
\left(\begin{array}{c} \lambda_{k+1,0}\\\lambda_{k+1,1} \end{array}\right)
= A \left(\begin{array}{c} \lambda_{k,0}\\\lambda_{k,1} \end{array}\right)
\end{equation*}
where 
\begin{equation*}
A=\left(\begin{array}{cc} 1&1/2\\1&1 \end{array}\right)
\end{equation*}
\begin{equation*}
\left(\begin{array}{c} \lambda_{k+1,0}\\\lambda_{k+1,1} \end{array}\right)
= A^k \left(\begin{array}{c} 1\\1 \end{array}\right)
\end{equation*}
Note that the maximum eigen value of $A$ is $1+1/\sqrt{2}$. Let $b=1+1/\sqrt{2}$
\begin{IEEEeqnarray*}{rCl}
\lambda_{k+1}&=&\lambda_{k+1,0}+\lambda_{k+1,1}\\
&\leq&c b^{k+1}
\end{IEEEeqnarray*}
for some $c\in\mathbb{R}^+$.
From equation \ref{maineqn}, we have 
\begin{IEEEeqnarray*}{rCl}
S_k&\leq&\frac{N^2}{k}\Omega_k\\
&\leq& \frac{N^2}{k} 2^{N-k}\lambda_k\\
&\leq& \frac{N^2}{k} 2^{N-k}c b^k
\end{IEEEeqnarray*}
To get a bound for $Q_1(k)$, we set this to be less than the total number of colorings $2^N$
\begin{IEEEeqnarray*}{rCl}
\frac{N^2}{k} 2^{N-k}c b^k&<&2^N\\
\Rightarrow N^2&<&\frac{k}{c}\left(\frac{2}{b}\right)^k\\
\Rightarrow N&<& \sqrt{\frac{k}{c}} \left(\sqrt{\frac{2}{b}}\right)^k \\
\Rightarrow N&<&c_0 g^k
\end{IEEEeqnarray*}
for some $c_0\in\mathbb{R}^+$ and $g=\sqrt{2/b}$. The value of $g$ can be calculated to be $1.08239..$.
\end{proof}
When rephrased in terms of probability, the results obtained in this chapter looked at the distribution of the number of monochromatic progressions in a random coloring. From calculating the expectation of the number of monochromatic progressions in a random coloring, we were able to deduce certain properties (more specifically, the existence of colorings avoiding monochromatic progressions). In the next chapter, we will look at more sophisticated techniques to analyse the distribution of random variables in general.

\chapter{Concentration Inequalities}

Concentration inequalities can be used to analyse how a random variable's distribution is concentrated. In this chapter, we will look at various concentration inequalities, their pros and cons, how to apply them and when to apply them. Alon and Spencer's book \cite{alon2011probabilistic} contains more of such inequalities and examples. This chapter arose out of reading Alon and Spencer's book and contains many interesting results, including, but not restricted to, the ones encountered when reading the book. The focus will be on application to problems. Given a problem, we would have to define an appropriate random variable that satisfies the hypothesis of some concentration inequality, then apply the inequality and finally analyse the result we get.

As each inequality will have different hypothesis on the random variable, it is not possible to have a single example which when attacked with different inequalities, will give us a complete understanding of all of them. So, it is necessary to give different examples for each inequality. Also note that the random variable we choose will depend on the problem and the inequality that we intend to use. Let us start with the simplest examples of Chebyshev and Chernoff bounds.
\section{Chebyshev and Chernoff Inequalities}\label{chebcher}
These are some of the simplest inequalities. They do not assume a lot on the random variable.
\subsection{Chebyshev Inequality}
\begin{theorem}
For any random variable $X$,
\begin{equation*}
Pr[|X-\mu|\geq \lambda\sigma]\leq \frac{1}{\lambda^2}
\end{equation*}
\end{theorem} 
\begin{proof}
\begin{IEEEeqnarray*}{rCl}
Var(X)&=&E[(X-\mu)^2]\\
&\geq& \lambda^2\sigma^2 Pr[|X-\mu|\geq \lambda\sigma]
\end{IEEEeqnarray*}
But $Var(X)=\sigma^2$.
\end{proof}

If we calculate the asymptotics when $X$ is a normal distribution, we get
\begin{equation*}
Pr[|X-\mu|\geq \lambda\sigma]=\sqrt{\frac{2}{\pi}}\frac{e^{-\lambda^2/2}}{\lambda}
\end{equation*}
This is significantly smaller than $1/\lambda^2$ (which we get from Chebychev inequality).
\begin{example}
Let us consider the following random variable
\begin{equation*}
X = \left\{
  \begin{array}{ll}
    0 & \text{with probability } 1-2p\\
    a & \text{with probability } p \\
    -a & \text{with probability } p
  \end{array}
\right.
\end{equation*}
$E[X]=0$, $Var[X]=E[X^2]=2pa^2$. Using Chebychev inequality, we get
\begin{equation*}
Pr[|X-\mu|\geq a]\leq \frac{\sigma^2}{a^2}=2p
\end{equation*}
But we know from the definition of the random variable that $Pr[|X|\geq a]=2p$. Hence, the Chebychev bound is tight for this example.
\end{example}
This inequality does not assume anything on the random variable. Hence, in cases where the random variable has more structure, we may not get tight bounds.
\subsection{Chernoff inequality}
Let $X=\sum_{i=1}^nX_i$, where the $X_i$'s are independent. I will outline the main idea in the proof of Chernoff inequality.
\begin{equation*}
Pr[X>(1+\delta)\mu]=Pr[e^{tX}>e^{(1+\delta)t\mu}]\leq \frac{E(e^{tX})}{e^{(1+\delta)t\mu}}
\end{equation*}
Note that
\begin{equation*}
E[e^{tX}]=E[e^{t\sum X_i}]=E[\prod_ie^{tX_i}]=\prod_iE[e^tX_i]
\end{equation*}
The last step is because we assumed that the $X_i$'s are independent. Given a problem, depending on the properties of $X_i$, we have to calculate $E[e^{tX_i}]$, plug it back and optimise over $t$ to get the best bound.

Let $X=\sum_{i=1}^nX_i$, where the $X_i$'s are independent. Each $X_i$ takes the value 1 with probability $p_i$ and 0 otherwise. Then, $E[X_i]=p_i$, $\mu=E[X]=\sum_{i=1}^np_i$. Applying the method outlined above, we get
\begin{IEEEeqnarray*}{rCl}
Pr[X>(1+\delta)\mu]&\leq& e^{-\frac{\delta^2\mu}{3}} \quad (0<\delta<1)\\
Pr[X<(1-\delta)\mu]&\leq& e^{-\frac{\delta^2\mu}{2}} \quad (0<\delta<1)
\end{IEEEeqnarray*}
\begin{example}
Consider $n$ coin flips of an unbiased coin. Let $X$ be the random variable giving the number of heads obtained. Then, $X=\sum_{i=1}^nX_i$. $X_i$'s are indicator random variables. They give 1 if the $i^{\text{th}}$ flip is head and 0 if it is a tail. $X_i$'s are independent, $E[X_i]=1/2$ and $\mu=E[X]=n/2$.
\begin{equation*}
Pr[X\geq \mu+\lambda]=Pr[X\geq \mu\left(1+\frac{\lambda}{\mu}\right)]\leq e^{-\left(\frac{\lambda}{\mu}\right)^2\frac{\mu}{3}}=e^{-\frac{\lambda^2}{3\mu}}
\end{equation*}
Now, let us choose $\lambda=O(\sqrt{n\log n})$. Then,
\begin{equation*}
Pr[X\geq \mu+\lambda]\leq \frac{1}{n^{2/3}}
\end{equation*}
\end{example}
Note that applying Chebychev inequality to the above example, we get, $Pr[X\geq \mu +\lambda]\leq \frac{\sigma^2}{\lambda^2}=\frac{1}{\log n}$. From this we see that when the random variable we consider has certain properties, it is best to apply inequalities that are specifically designed for their case. In the example, we had independence of $X_i$'s. Using the Chebychev inequality did not use this property. So, we got a weak bound. When we used the Chernoff inequality, we used this information and hence got a much stronger bound.
\section{Azuma's Inequality}
We will start this section by proving Azuma's inequality and then look at its applications. Let us recall some definitions.

\begin{definition}
A martingale is a sequence of random variables $X_0,X_1,\ldots,X_m$ such that $E[X_{i+1}|X_i,X_{i-1},\ldots X_0]=X_i$ for $0<i<m$.
\end{definition}
\subsection{Basics}
\begin{theorem}[Azuma]
Let $X_0,X_1,\ldots,X_n$ be a martingale with $|X_{i+1}-X_i|\leq 1$ for $0\leq i<m$. Then,
\begin{equation*}
Pr[X_m > \lambda \sqrt{m}]<e^{-\lambda^2/2}
\end{equation*}
for $\lambda\in\mathbb{R}^+$.
\end{theorem}
\begin{proof}
Let $Y_i=X_i-X_{i-1}$ for $0<i\leq n$. Then $|Y_i|\leq 1$ and $E[Y_i|X_{i-1},\ldots X_0]=0$.
If we define $h(x)=\cosh \lambda +x \sinh \lambda$, we get,
\begin{equation*}
E[e^{\alpha Y_i}]\leq E[h(Y_i)]=h(E[Y_i])=h(0)=\cosh \lambda
\end{equation*}
In the previous analysis, we can look at $E[e^{\alpha Y_i}|X_{i-1},X_{i-2},\ldots,X_0]$ instead of $E[e^{\alpha Y_i}]$, and get the following
\begin{equation*}
E[e^{\alpha Y_i}|X_{i-1},X_{i-2},\ldots,X_0] \leq \cosh \lambda\leq e^{\alpha^2/2}
\end{equation*}
Note that $\cosh \lambda\leq e^{\alpha^2/2}$ follows directly from comparing their power series. Now,
\begin{IEEEeqnarray*}{rCl}
E[e^{\alpha X_m}]&=&E\left[\prod_{i=1}^me^{\alpha Y_i}\right]\\
&=&E\left[\left(\prod_{i=1}^{m-1}e^{\alpha Y_i}\right) E[e^{\alpha Y_m}|X_{m-1},X_{i-2},\ldots,X_0] \right]\\
&\leq&E\left[\left(\prod_{i=1}^{m-1}e^{\alpha Y_i}\right) \right] e^{\alpha^2/2}\\
&\leq&e^{m\alpha^2/2}
\end{IEEEeqnarray*}
Now, if we apply Markov inequality to $Pr[e^{\alpha X_m}>e^{\alpha\lambda \sqrt{m}}]$, and optimise over $\alpha$, we get the desired result.
\begin{IEEEeqnarray*}{rCl}
Pr[X_m>\lambda\sqrt{m}]&=&Pr[e^{\alpha X_m}>e^{\alpha\lambda \sqrt{m}}]\\
&<&E[e^{\alpha X_m}]e^{-\alpha\lambda\sqrt{m}}\\
&\leq&e^{\alpha^2m/2-\alpha\lambda\sqrt{m}}\\
&=&e^{-\lambda^2/2}
\end{IEEEeqnarray*}
\end{proof}
\begin{remark}
Under the same hypothesis of the previous theorem, we can also get that
\begin{equation*}
Pr[X_m < -\lambda \sqrt{m}]<e^{-\lambda^2/2}
\end{equation*}
In general, for a martingale $c=X_0,X_1,X_2,\ldots,X_m$ with $|X_{i+1}-X_i|\leq 1$ for $0\leq i<m$, we get
\begin{equation*}
Pr[|X_m-c| > \lambda \sqrt{m}]<2e^{-\lambda^2/2}
\end{equation*}
for $\lambda\in\mathbb{R}^+$.
\end{remark}

We will be looking at applications of Azuma's inequality to random graphs. A random graph $G(n,p)$ is a graph on $n$ labelled vertices obtained by selecting each pair of vertices to be an edge with probability $p$ randomly and independently. We have to define a martingale on a random graph. We will define two different martingales. Depending on the problem, we will select one of them.

\begin{definition}
Let $G(n,p)$ be the underlying probability space. Label the unordered pairs of vertices $\{i,j\}$ by $e_1,e_2,\ldots,e_m$ where $m={n\choose 2}$. These represent the possible edges. We define a martingale $X_0,X_1,\ldots, X_m$ the following way. Let $f$ be any graph theoretic function. For $H\in G(n,p)$, $X_m(H)=f(H)$, $X_0(H)=E[f(G)]$ and in general,
\begin{equation*}
X_i(H)=E[f(G)|e_j\in G \Leftrightarrow e_j\in H, 1\leq j\leq i]
\end{equation*}
\end{definition}
To find $X_i(H)$, we expose the first $i$ edges of $H$ (the remaining edges are considered to be randomly chosen with probability $p$), and compute the expectation of $f(G)$ with this information. So, each successive random variable has more information about $f(G)$. The fact that the remaining edges are considered random and $X_i(H)$ gives the expectation of $f(G)$ under this partial information automatically makes the $X_i$'s a martingale. $f$ is said to satisfy the edge Lipschitz condition if whenever $H$ and $H'$ differ in only one edge then, $|f(H)-f(H')|\leq1$. The martingale that arises from such an $f$ satisfies $|X_{i-1}-X_i|\leq1$.

\begin{definition}
Let $G(n,p)$ be the underlying probability space. We define a martingale $X_1,\ldots, X_n$ the following way. Let $f$ be any graph theoretic function. For $H\in G(n,p)$,
\begin{equation*}
X_i(H)=E[f(G)| \forall x,y\leq i, \{x,y\}\in G\Leftrightarrow \{x,y\}\in H]
\end{equation*}
\end{definition}
To find $X_i(H)$, we expose the first $i$ vertices and all their internal edge and takes the conditional expectation of $f(G)$ with this information. Note $X_1(H)=E[f(G)]$ and $X_n(H)=f(H)$. $f$ is said to satisfy the vertex Lipschitz condition if whenever $H$ and $H'$ differ at only one vertex, $|f(H)-f(H')|\leq1$. The martingale that arises from such an $f$ satisfies $|X_{i-1}-X_i|\leq1$.

\subsection{Applications}

Before we start, we need the following definition.

\begin{definition}
The chromatic number of a graph $G$ (denoted $\chi(G)$) is the smallest number of colors needed to color the vertices of $G$ so that no two adjacent vertices have the same color.
\end{definition}

The study of the chromatic number of random graphs was pioneered by Erd\H{o}s and R\'enyi \cite{erdos1960erg}. In theorem \ref{shamirthm}, once we define the martingale, the result will follow immediately from Azuma's inequality. This result (proved in \cite{shamir}) was the first to use the martingale approach in probabilistic method. This shows that the chromatic number of a random graph is concentrated in an interval of length $O(\sqrt{n})$. 

Bollob{\'a}s \cite{bollobas1988chromatic} proved for any $p\in (0,1)$, almost surely,
\begin{equation*}
\chi(G(n,p))=\frac{n}{2\log n}\log\left(\frac{1}{1-p}\right)(1+o(1)).
\end{equation*}
This result again used Azuma's inequality to first bound the size of maximum clique of $G$. This is a more complicated application of Azuma's inequality and its proof is outlined in theorem \ref{bollobasthm}.

Going back to concentration of $\chi(G)$, in theorem \ref{4conc}, we prove that $\chi(G(n,p))$ is concentrated in an interval of length 4 when $p=n^{-\alpha}$ and $\alpha>5/6$. When $p<n^{-1/6-\epsilon}$, {\L}uczak \cite{luczak1991note} showed that the value of $\chi(G(n,p))$ is concentrated in an interval of length 2. Alon and Krivelevich extended this result for $p<n^{-1/2-\epsilon}$ in \cite{alon1997concentration}. These results only give a range of interval where $\chi(G)$ will mostly lie. They do not say where the interval is. Recently, this question was answered by Achlioptas and Naor \cite{zbMATH05033810}. Now, let us start with our examples.

\begin{theorem}
Let $G\sim G(n,p)$ for some fixed $n\in \mathbb{Z}^+$, $p\in[0,1]$. Then, \label{shamirthm}
\begin{equation*}
Pr[|\chi(G)-E[\chi(G)]|>\lambda \sqrt{n-1}]< 2 e^{-\lambda^2/2}
\end{equation*}
\end{theorem}
\begin{proof}
Let $X_1,X_2,\ldots,X_n$ be the vertex exposure martingale on $G(n,p)$ with $f(G)=\chi(G)$. This function satisfies the vertex Lipschitz condition. Hence, applying Azuma's inequality, we get the result directly.
\end{proof}

Now, let us prove a bound for $\omega (G)$ (size of the maximum clique of $G$). This result was first proved by B\'ela Bollob\'as. A clique $C$ is a collection of vertices of $G$ where for every $i,j\in C$, there is an edge between $i$ and $j$. First let us fix $k\in\mathbb{Z}^+$. We will consider $Y=Y(H)$ to be the maximal size of family of edge disjoint cliques of size $k$ in $H$.
\begin{theorem} \label{bollobasthm}
Let $n\in \mathbb{Z}^+$ be the number of vertices, $k\sim 2\log_2n$ and $G\sim G(n,p)$. Then,
\begin{equation*}
Pr[\omega(G)<k]<e^{(c+o(1))(n^2/k^8)}
\end{equation*}
\end{theorem}
\begin{proof}
Let $Y_0,Y_1,\ldots,Y_m$ where $m={n\choose 2}$ be the edge exposure martingale on $G(n,1/2)$ with the function $Y$. Note that $Y$ satisfies the edge Lipschitz condition. We are going to apply Azuma's inequality and get a bound on the probability that $Y_m=0$. Note that $Y_m=0$ is the same as saying that $\omega(G)<k$. Hence,
\begin{IEEEeqnarray*}{rCl}
Pr[\omega(G)<k]&=& Pr[Y_m=0]\\
&\leq&Pr[Y_m-E[Y]\leq -E[Y]]\\
&\leq& e^{-E[Y]^2/2{n\choose 2}}\\
&\leq&e^{(c+o(1))(n^2/k^8)}
\end{IEEEeqnarray*}
The last step uses $E[Y]\geq (1+o(1))\frac{n^2}{2k^4}$. This can be proved using the probabilistic method under the assumptions of this theorem.
\end{proof}
Next, we will give a more difficult application of this inequality. We will need the following lemma. I will only outline the proof of the lemma.
\begin{lemma}
Let $\alpha,c$ be fixed, $\alpha>5/6$ and $p=n^{-\alpha}$. Then, the probability that $c\sqrt{n}$ vertices of $G=G(n,p)$ is 3-colorable tends to 1 as $n$ tends to infinity.
\end{lemma}
\begin{proof}
Let $T$ be the minimal set that is not 3-colorable. Each vertex of $T$ must have internal degree at least 3. Then, $T$ has at least $3t/2$ edges. Now, the probability that this happens for some $T$ with at most $c\sqrt{n}$ vertices is bounded above by
\begin{equation*}
\sum_{t=4}^{c\sqrt{n}}{n\choose t}{{t\choose 2}\choose 3t/2}p^{3t/2}
\end{equation*}
With careful bounding, we can show that this is $o(1)$. Hence, $Pr[|T|>c\sqrt{n}]$ goes to 0 as $n$ goes to infinity.
\end{proof}
Now, let us prove the following theorem.
\begin{theorem}
Let $p=n^{-\alpha}$, where $\alpha>5/6$ and $G=G(n,p)$. Then, there exists $u=u(n,p)$ such that almost always \label{4conc}
\begin{equation*}
u\leq\chi(G)\leq u+3
\end{equation*}
\end{theorem}
\begin{proof}
Let $\epsilon>0$ be arbitrary and $u=u(n,p,\epsilon)$ be the least integer such that $Pr[\chi(G)\leq u]>\epsilon$. Let $S(G)$ be the minimal set of vertices $S$ such that $G-S(G)$ is $u$-colorable and $Y(G)$ be the size of $S(G)$. Note that by this choice of $u$, $Pr[Y=0]>\epsilon$. For any subset $B$ of the vertices of $G$, let $A_B$ be the event that the set $B$ can not be 3-colored. If we show that $Pr[Y\geq c\sqrt{n}]\leq \epsilon$, we also know when $Y<c\sqrt{n}$, these vertices can be 3-colored almost always. Hence, we can get that
\begin{IEEEeqnarray*}{rCl}
Pr[\chi(G)\leq u+3]&=& 1-Pr[A_{S(G)}]\\
&\geq&1-Pr[Y\geq c\sqrt{n}]-Pr[A_{S(G)}|Y<c\sqrt{n}]\\
&\geq&1-2\epsilon
\end{IEEEeqnarray*}
Further $Pr[\chi(G)<u]<\epsilon$. Hence, $Pr[u\leq \chi(G)\leq u+3]\geq 1-3\epsilon$. Now, we have completed the proof except for proving $Pr[Y\geq c\sqrt{n}]\leq \epsilon$. This is the part that uses Azuma's inequality. Consider the vertex exposure martingale defined by Y on $G(n,p)$. This function satisfies the vertex Lipschitz condition. Hence,
\begin{IEEEeqnarray}{rCl}
Pr[Y\leq \mu-\lambda\sqrt{n-1}]&<&e^{-\lambda^2/2}\label{firsteqn}\\
Pr[Y\geq \mu+\lambda\sqrt{n-1}]&<&e^{-\lambda^2/2}\label{secondeqn}
\end{IEEEeqnarray}
Let $\lambda$ satisfy $e^{-\lambda^2/2}=\epsilon$. Then, $Pr[Y=0]>\epsilon$ and $Pr[Y\leq \mu-\lambda\sqrt{n-1}]<\epsilon$ implies that $\mu-\lambda\sqrt{n-1}<0\Rightarrow\mu< \lambda\sqrt{n-1}$. Substituting this into \ref{secondeqn}, we get that
\begin{equation*}
Pr[Y\geq 2\lambda\sqrt{n-1}]\leq Pr[Y\geq \mu +\lambda \sqrt{n-1}]\leq \epsilon
\end{equation*}
\end{proof}
\section{Talagrand Inequality}

Let $\Omega_1,\Omega_2,\ldots,\Omega_n$ be probability spaces and $\Omega$ denote the product space. $\{X_i\}$ is a family of independent random  variables with each $X_i$ taking values in $\Omega_i$. If $x_1\in X_1, x_2\in X_2,\ldots x_n\in X_n$, then $x=(x_1,x_2,\ldots,x_n)\in\Omega$. Let us define a notion of distance in $\Omega$ as follows. Let $A\subseteq \Omega$. Then, the distance from $x$ to $y\in A$ is defined as 
\begin{equation*}
d(x,y)=\sup_{|\alpha|=1}\sum_{x_i\neq y_i}\alpha_i
\end{equation*}
where $\alpha$ varies over all unit vectors in $\mathbb{R}^n$. We also define
\begin{equation*}
d(x,A)=\inf_{y\in A} d(x,y)
\end{equation*}
One more definition that we will need before stating Talagrand's theorem is \\
$A_t=\{x\in \Omega\quad:\quad d(x,A)\leq t\}$
Now, we are ready to state Talagrand's inequality. The proof uses induction on the dimension of $\Omega$ and can be found in \cite{alon2011probabilistic}.
\begin{theorem}
Let $\{X_i\}_{i=1}^n$ be independent random variables such that $X_i$ takes values from $\Omega_i$ and $A\subseteq\Omega$. Then, for any $t\geq0$
\begin{equation*}
Pr[A](1-Pr[A_t])\leq e^{-t^2/4}
\end{equation*}
\end{theorem}
I will not restrict to applications in graph theory alone. Before going to applications, I will state a corollary of the theorem that will be very useful. The corollary needs the following definition. Let $X=h(x)$ be a random variable that is Lipschitz. $h$ is said to be Lipschitz if $h:\Omega \rightarrow \mathbb{R}$ and $|h(x)-h(y)|\leq1$ whenever $x,y$ differ in at most one coordinate. Generalising this idea, we say a function $h$ is $k$-Lipschitz if $|h(x)-h(y)|\leq k$ whenever $x,y$ differ in only one coordinate. In applications, we would also demand that $h$ satisfies the following property.
\begin{definition}
Let $f:\mathbb{N}\rightarrow \mathbb{N}$. $h$ is $f$-certifiable if whenever $h(x)\geq s$, there exists $I\subseteq \{1,2,\ldots,n\}$ with $|I|\leq f(s)$ so that all $y\in \Omega$ that agree with $x$ on the coordinates of $I$ have $h(y)\geq s$.
\end{definition}
\begin{corollary}
Let $X$ be a random variable in a space $\Omega$ that is itself a product space. Suppose $X$ arises out of a function $h: \Omega\rightarrow \mathbb{R}$ that is $k$-Lipschitz and $f$-certifiable. Then, for all $b,t$,
\begin{equation*}
Pr[X\leq b-tk\sqrt{f(b)}]Pr[X\geq b]\leq e^{-t^2/4}
\end{equation*}
\end{corollary}
If $b$ is chosen to be the median ($m$), we get a bound for $X$ falling below $t\sqrt{f(b)}$ of the median. If we want a bound for $X$ going above $t\sqrt{f(b)}$ of the median, we should choose $m=b-t\sqrt{f(b)}$.
I start with a simple example.
\begin{example}
Let $\Omega$ be the probability space that is uniformly distributed in $[0,1]^n$. $\Omega$ can be considered as the product space of $\Omega_i$. Where each $\Omega_i$ gets values from $[0,1]$ uniformly and independently. Let $X=h(x)$ be the random variable length of the longest increasing subsequence in $x=(x_1,x_2,\ldots,x_n)$.

Note that $h(x)$ is 1-Lipschitz and $f$-certifiable with $f(s)=s$ since if $x$ has an increasing subsequence of length $s$, these $s$ coordinates will certify that $X\geq s$. Hence, from Talagrand's inequality, we get that
\begin{IEEEeqnarray*}{rCl}
Pr[X<m-t\sqrt{m}]&\leq& 2e^{-t^2/4}\\
Pr[X>m+t\sqrt{m}]&\leq& 2e^{-t^2/4}
\end{IEEEeqnarray*}
We can also show that $m=\Theta(n^{1/2})$. Hence, $|X-m|<n^{1/4}$ almost surely.
\end{example}

\begin{example}
Let us find a bound for $Pr[\omega(G)<k]$ using Talagrand's inequality. Recall that $k\sim 2\log_2n$ and $G \sim G(n,1/2)$. The bound we got using Azuma's inequality was
\begin{equation*}
Pr[\omega(G)<k]<e^{(c+o(1))(n^2/k^8)}
\end{equation*}
We used Azuma's inequality on the random variable $Y$ which was defined to be the maximal number of edge disjoint $k$ cliques in $G$. We will use the same random variable and apply Talagrand's inequality. Before that note that the probability we require is the same as the probability that $Y\leq 0$.
\begin{equation*}
Pr[\omega(G)<k]=Pr(Y\leq 0)
\end{equation*}

We know that $Y$ is tightly concentrated about $n^2/k^4$. Hence, $m\sim n^2/k^4$. $Y$ is 1-Lipschitz and $f$-certifiable with $f(s)={k\choose 2}s$. Hence, using Talagrand's inequality, we get
\begin{equation*}
Pr[Y\leq m-tm^{1/2}{k\choose 2}^{1/2}]Pr[Y\geq m]<e^{-t^2/4}
\end{equation*}
We choose $t$ such that $m-tm^{1/2}{k\choose 2}^{1/2}=0$. Solving for $t$, we get $t=\Theta(n/\ln^3n)$. Then
\begin{equation*}
Pr[\omega(G)<k]=Pr[Y\leq 0]<2e^{-t^2/4}<e^{-\Omega(n^2/\ln^6n)}
\end{equation*}
\end{example}
\section{Janson's Inequality}
Let $B_i$'s be rare events in some probability space. We want to show that $Pr[\wedge_i \overline{B_i}]$ can be made very small. If the events are independent, it follows directly that
\begin{equation*}
Pr[\wedge_i \overline{B_i}]=\prod_iPr[\overline{B_i}]=M \text{ (say) }
\end{equation*}
But if there is dependency between the $B_i$'s, we can still use Janson's inequality to show that $Pr[\wedge_i \overline{B_i}]$ lies close to $M$. How close it is will depend on how rare the events are and the significance of the dependencies.

Let us make all this formal. Let $\Omega$ be a finite set. Let $R$ be a random subset of $\Omega$ given by $Pr[r\in R]=p_r$. The elements are picked randomly and independently. Let $\{A_i\}_{i\in I}$ be subsets of $\Omega$ where $I$ is a finite index set. Let $B_i$ be the event that all $r\in A_i$ were picked and $X_i$ be the indicator random variable for $B_i$. $X=\sum_iX_i$ gives the number of $A_i\subseteq R$. For $i\neq j$, if $A_i \cap A_j\neq \phi $, the events $B_i$ and $B_j$ are dependent. Further, note that the occurrence of event $B_i$ will positively influence the occurrence of the event $B_j$. This is because if we know that the event $B_i$ occurred, all the elements of $A_i$ were picked. This would mean that some of the elements of $A_j$ were picked. So, with this information, the occurrence of $B_j$ has more probability. Similarly, when we know that the event $B_i$ did not occur, the probability that the event $B_j$ occurs goes down. This can be written as
\begin{equation}
Pr\left[B_i| \bigwedge_{j\in J}\overline{B_j}\right]\leq Pr[B_i] \label{correlation1}
\end{equation}
for all sets $J\subset I$ and $i\notin J$. In general, we can say that
\begin{equation}
Pr\left[B_i| B_k\wedge \bigwedge_{j\in J}\overline{B_j}\right]\leq Pr[B_i|B_k] \label{correlation2}
\end{equation}
for all sets $J\subset I$ and $i,k\notin J$. These correlation inequalities will play a critical role in proving Janson's Inequality. 
\subsection{Proof}
Before starting the proof, we require the following definitions.
\begin{IEEEeqnarray*}{rCl}
\Delta&=&\sum_{i\sim j}Pr[B_i\wedge B_j]\\
M&=&\prod_{i\in I}Pr[\overline{B_i}]\\
\mu&=&E[X]=\sum_{i\in I} Pr[B_i]
\end{IEEEeqnarray*}
\begin{theorem}
Let $\{B_i\}_{i\in I}$ be events as defined above with $Pr[B_i]\leq \epsilon$. Then,
\begin{equation*}
M\leq Pr[\wedge_{i\in I}\overline{B_i}]\leq Me^{\frac{\Delta}{2(1-\epsilon)}}
\end{equation*}
\end{theorem}
From the theorem, we can say that $Pr[\wedge_{i\in I}\overline{B_i}]\sim M$ when $e^{\frac{\Delta}{2(1-\epsilon)}}\rightarrow 1$. This happens when $\epsilon$ and $\Delta$ are small.
\begin{proof}
Let $I=\{1,2,3,\ldots,m\}$. Then,
\begin{equation*}
Pr[\wedge_{i\in I}\overline{B_i}] = \prod_{i=1}^m Pr[\overline{B_i}| \wedge_{1\leq j<i}]\geq \prod_{i=1}^mPr[\overline{B_i}]
\end{equation*}
The last step follows from \ref{correlation2}. Now for the upper bound, the idea is to find an upper bound for $Pr[\overline{B_i}|\wedge_{1\leq j<i}\overline{B_j}]$ and then take the product over all $i$. Then, the theorem will follow from the following
\begin{equation}
Pr[\wedge_{i\in I}\overline{B_i}]=\prod_{i-1}^mPr[\overline{B_i}|\wedge_{1\leq j<i}\overline{B_j}] \label{jansoneqn}
\end{equation}

For every $i$, renumber the first $i-1$ elements so that $i\sim j$ for $1\leq j\leq d$ and not so for the remaining elements. Note that $d$ will depend on $i$. Let $A=B_i, B=\overline{B_1}\wedge\cdots\wedge \overline{B_d}$ and $C=\overline{B_{d+1}}\wedge\cdots\wedge \overline{B_{i-1}}$. 

\begin{IEEEeqnarray*}{rCl}
Pr[B_i|\wedge_{1\leq j<i}\overline{B_j}]&=&Pr[A|B\wedge C]\\
&\geq & Pr[A\wedge B|C]\\
&=& Pr[A|C] Pr[B|A\wedge C]\\
&=& Pr[B_i]Pr[B|A\wedge C]
\end{IEEEeqnarray*}
Now, we have to bound $Pr[B|A\wedge C]$.
\begin{IEEEeqnarray*}{rCl}
Pr[B|A\wedge C] &\geq& 1-\sum_{j=1}^d Pr[B_j|B_i\wedge C]\\
&\geq & 1-\sum_{j=1}^d Pr[B_j|B_i]
\end{IEEEeqnarray*}
Hence, we get,
\begin{IEEEeqnarray*}{rCl}
Pr[\overline{B_i}|\wedge_{1\leq j<i}\overline{B_j}]&\leq& Pr[\overline{B_i}]+ \sum_{j=1}^d Pr[B_j\wedge B_i]\\
&\leq& Pr[\overline{B_i}]\left(1+\frac{1}{1-\epsilon} \sum_{j=1}^d Pr[B_j\wedge B_i]\right)\\
&\leq& Pr[\overline{B_i}]\exp\left( \frac{1}{1-\epsilon} \sum_{j=1}^d Pr[B_j\wedge B_i] \right)
\end{IEEEeqnarray*}
Substituting this into \ref{jansoneqn}, we get the desired result.
\end{proof}

\subsection{Applications}
Let us start with a simple example.
\begin{theorem}
Let $G\sim G(n,p)$. Let $c$ be some constant and $p=c/n$. Let $A$ be the event that there is no triangle in $G$. Then, asymptotically,
\begin{equation*}
Pr[A]\rightarrow e^{-c^3/6}
\end{equation*}
\end{theorem}
\begin{proof}
Let $G\sim G(n,p)$. Let $\{S_i\}$ be a collection all possible triangles in $G$. There are ${n\choose 3}$ ways to choose a possible triangle. Hence, $i$ varies from $1$ to ${n\choose 3}$. Let $B_i$ be the event that $S_i$ is a triangle in $G$. These are the rare (bad) events. If none of these events happen, we can say that $A$ has happened.
\begin{equation*}
Pr[B_i]=p^3=\frac{c^3}{n^3}=o(1)
\end{equation*}
$p^3$ can be chosen to be $\epsilon$. Further,
\begin{equation*}
M=\prod_iPr[\overline{B_i}]=(1-p^3)^{n\choose 3}\rightarrow e^{-\frac{c^3}{6}}
\end{equation*}
And,
\begin{IEEEeqnarray*}{rCl}
\Delta&=& \sum_{i\sim j} Pr[B_i\wedge B_j]\\
&=& \sum_{i\sim j} p^5\\
&=&{n\choose 3}{3\choose2}{n-3 \choose 1}p^5\\
&=&O(n^5p^5)\\
&=&O\left(n^4\frac{c^5}{n^5}\right)=o(1)
\end{IEEEeqnarray*}
Hence, applying Janson's inequality, we get that
\begin{equation*}
M\leq Pr[A]\leq Me^{o(1)}
\end{equation*}
Asymptotically, $Pr[A] \rightarrow e^{-c^3/6}$.
\end{proof}
Now, let us look at a more complicated example.
\begin{theorem}
Let $G\sim G(n,p)$. Let $B$ be the event that there exists a path of length 3 between any pair of vertices of $G$. Then, $Pr[B]\rightarrow 1$ asymptotically when $p=\left(\frac{c\ln n}{n^2}\right)^{1/3}$ and $c\geq 2$.
\end{theorem}
\begin{proof}
Let $u,v$ be vertices in $G$ and $B_{u,v}$ be the event that there does not exist a path of length 3 between $u$ and $v$. Then,
\begin{IEEEeqnarray*}{rCl}
Pr[\overline{B}]&=&Pr[\cup_{u,v}B_{u,v}]\\
&\leq& \sum_{u,v}Pr[B_{u,v}]\\
&=&O(n^2)Pr[B_{u,v}]
\end{IEEEeqnarray*}
Now, our aim is to use Janson's inequality to prove $Pr[B_{u,v}]=o(n^{-2})$. When considering $B_{u,v}$, the bad events will be $A_{w_1,w_2}$. $A_{w_1,w_2}$ is the event that there is a path of length 3 in $G$ from $u$ to $v$ that goes from $u$ to $w_1$, $w_1$ to $w_2$ and $w_2$ to $v$ ($w_1$ and $w_2$ are vertices in $G$). In other words, $A_{w_1,w_2}$ is the event that the following edges are present in $G$: $(u,w_1)$, $(w_1,w_2)$ and $(w_2,v)$. Clearly,
\begin{equation*}
Pr[B_{u,v}]=Pr[\wedge_{w_1,w_2}\overline{A_{w_1,w_2}}]
\end{equation*}
As we are only interested in the asymptotic behaviour, we will calculate the quantities necessary for applying Janson's inequality asymptotically.
\begin{equation*}
Pr[A_{w_1,w_2}]=p^3=\frac{c\ln n}{n^2}\rightarrow 0
\end{equation*} 
This is independent of the choice of vertices $w_1$ and $w_2$. We can choose this to be our $\epsilon$. Next,
\begin{IEEEeqnarray*}{rCl}
M&=&\prod_{w_1,w_2}Pr[\overline{A_{w_1,w_2}}]\\
&=&(1-p^3)^{(n-2)(n-3)}\\
&\leq& e^{-p^3n^2}\\
&=&n^{-c}=o(n^{-2})
\end{IEEEeqnarray*}
Similarly, for calculating $\Delta$, we should consider the three possible ways two paths can intersect and sum their probability. The contribution of the leading term will be
\begin{IEEEeqnarray*}{rCl}
\Delta&=&\sum_{(w_1,w_2)\sim (w_1',w_2')}Pr[A_{w_1,w_2}\wedge A_{w_1',w_2'}]\\
&=&O(n^3)p^5\\
&=&O\left( \frac{(c\ln n)^{5/3}}{n^{1/3}}\right)=o(1)
\end{IEEEeqnarray*}
Hence,
\begin{IEEEeqnarray*}{rCl}
Pr[B_{u,v}]=Pr[\wedge_{w_1,w_2}\overline{A_{w_1,w_2}}]=o(n^{-2})
\end{IEEEeqnarray*}
\end{proof}

\section{Summary}

In this chapter, we started with simple inequalities like Chebychev inequality and Chernoff inequality and then moved to more complicated ones like Azuma's inequality, Talagrand's inequality and Janson's inequality. Note that this is not an exhaustive list of all concentration inequalities. 

In section \ref{chebcher}, I discussed how all the information present in the random variable has to be used in order to get the best results. In the example I considered, when the independence information was not used, we got much weaker results. Given a problem, the first and most important thing to do (in order to apply concentration inequality) is to define the appropriate random variable. This has to be done with a concentration inequality in mind because the random variable has to satisfy the hypothesis of the inequality. This should also be done in such a way that all the information in the problem will be used when applying the inequality. This is the toughest part when using such inequalities. After defining the appropriate random variable, it is a direct application of inequality to get the result. So, if we ignore some information in the problem, we will end up getting much weaker results. 

Applying these inequalities along with the probabilistic method to problems similar to the ones considered in the first three chapters of this thesis will give us more information about the distribution of the random variable considered.

\bibliographystyle{plain}
\bibliography{report}

\end{document}